\documentstyle[12pt]{article} \textwidth 6.1 in \textheight 9.5 in

\title{\bf Asymptotic behavior of positive solutions of semilinear elliptic equations in $R^{n}$
\\
  \thanks{This work was supported  by the XiaoXiang Funds of Hunan Normal
 University, China and the Natural Science Foundation of
China (10671064) and   Education Foundation  of Hunan Province
(06c516).} }
\author{Baishun Lai,  Shuqing Zhou, qing Luo\\
} \baselineskip 0.2in
\date{}
\begin{document}
\maketitle
\begin{center}
\begin{minipage}{130mm}
{\small {\bf Abstract}\vskip0.1in

\ \ \  We will investigate the asymptotic behavior of positive
solutions of the elliptic equation
\begin{equation}
\Delta u+|x|^{l_{1}}u^{p}+|x|^{l_{2}}u^{q}=0\ \mbox{in}\ \ R^{n}.\
\end{equation}
We establish that for $n\geq 3$ and $q>p>1$, any  positive radial
solution of (0.1) has the following property:
$\lim_{r\to\infty}r^{\frac{2+l_{1}}{p-1}}u$ and
$\lim_{r\to0}r^{\frac{2+l_{2}}{q-1}}u$  always exist if
$\frac{n+l_{1}}{n-2}<p<q, \ \ p\neq\frac{n+2+2l_{1}}{n-2},\ \ q
\neq\frac{n+2+2l_{2}}{n-2}.$ In addition, we prove that the singular
solution of (0.1) is unique under a certain condition.}
\end{minipage}
\end{center}
\baselineskip 0.2in \vskip 0.2in {\bf 1. Introduction} %
\setcounter{section}{1} \setcounter{equation}{0}\vskip 0.2in In this
paper we will study the asymptotic behavior of positive solutions of
the following equation
\begin{equation}
\Delta u + K_{1}(|x|)u^{p}+K_{2}(|x|)u^{q}=0 \ \ \ \mbox{in} \ \ \
R^{n},\ \ \ n\geq3,
\end{equation}
and in particular, of positive radial solutions of
\begin{equation}
\Delta u +|x|^{l_{1}}u^{p}+|x|^{l_{2}}u^{q}=0\ \ \ \mbox{in} \ \ \
R^{n},\ \ \ n\geq3,
\end{equation}
where $-2<l_{2}< l_{1}\leq0, 1<p< q$ and
$\Delta=\Sigma_{1}^{n}\frac{\partial^{2}}{\partial x_{i}^{2}}$ is
the Laplace operator, and $K_{1}, K_{2}$ are a locally H\"{o}lder
continuous in $R^{n} \setminus\{0\}$. By an entire solution of
(1.1), we mean a positive weak solution of (1.1) in $R^{n}$
satisfying (1.1) pointwise in $R^{n}\setminus \{0\}$.\vskip0.1in

When $K_{1}(|x|)=K_{2}(|x|), p=q$, then (1.2) reduces to, by scaling
\begin{equation}
\Delta u + K(|x|)u^{p}=0\ \ \ \mbox{in} \ \ \ R^{n},\ \ \ n\geq3.
\end{equation}
Equation (1.3) has their roots from many mathematical and physical
fields, e.g., the well-known scalar curvature equation in the study
of Riemannian geometry, the scalar field equation for the standing
wave of nonlinear Schr\"{o}dinger and the Klein-G\"{o}rder
equations, the Matukuma equation describing the dynamics of globular
cluster of star, etc. We refer the interested readers to [16, 18,
20, 21, 23] and the references therein. There have been many works
devoted to studying the existence of positive solutions of (1.3) in
$R^{n}$ after the first contribution by Ni [20] in 1982, see [2, 6,
7, 8, 11, 12, 23 ] and the references therein. One of the features
of the equation is that (1.3) can posses infinitely many solutions
as long as the exponent $p$ and the dimension $n$ are large enough.
Recent studies in [1, 3, 4, 9, 10 ] paid special attention to this
phenomenon. The purpose of this paper is to study the asymptotic
behavior of positive entire solutions and to present the uniqueness
result. In such perspective, we review related works as follows.  In
the fast decay case, i.e. $|K|\leq Cr^{l},
 l<-2$, Ni [20] showed that equation (1.3) possess infinitely many
positive solutions which are bounded from below by positive
constants.  Li and Ni [18] showed that, for positive solution of
(1.3), the limit $u_{\infty}=\lim_{x\to\infty}u(x)$ always exists.
Furthermore, if $u_{\infty}=0$, then, for any $\varepsilon>0$,
\begin{eqnarray*}
u(x)\leq \left\{
\begin{array}{ll}
C|x|^{2-n} & \ \ \ \mbox{if} \ \ \ p>\frac{n+l}{n-2},\\\\
C_{\varepsilon}|x|^{[(1-\varepsilon)(l+2)]/1-p}& \ \ \ \mbox{if} \ \
\ p\leq\frac{n+l}{n-2},
\end{array}
\right.
\end{eqnarray*}
where $C_{\varepsilon}$ is a constant depending on $\varepsilon$,
and if $u_{\infty}>0$, then
\begin{eqnarray*} |u(x)-u_{\infty}|\leq \left\{
\begin{array}{ll}
C|x|^{2-n} & \ \ \ \mbox{if} \ \ \ l<-n,\\\\
C|x|^{2-n}\log |x| & \ \ \ \mbox{if} \ \ \ l=-n,\\\\
C|x|^{2+l} & \ \ \ \mbox{if} \ \ \ -n<l<-2,
\end{array}
\right.
\end{eqnarray*}
at $\infty$.

For the slow decay case, i.e. $K(r)\geq Cr^{l},$ for some $l>-2$ and
$r$ large enough. And also $K(r)$ satisfy:

(K.1)\ \ $K(r)>0$ in $r>0$ and
$\lim_{r\to\infty}r^{-l}K(r)=k_{\infty}>0$,

(K.2)\ \ $K(r)$ is differentiable and
$[\frac{d}{dr}(r^{-l}K(r)]^{+}\in L^{1}, r>0$,

(K.2)\ \ $K(r)$ is differentiable and
$[\frac{d}{dr}(r^{-l}K(r)]^{-}\in L^{1}, r>0$.\vskip 0.2in
 Li [17]
gave an accurate description on the asymptotic behavior of positive
solutions of (1.3), which is stated as follows\vskip 0.1in

Theorem A (Li [ 17, Theorem 1]) \ Let $u$ be a positive radial
solution of (1.3). Assume that $K$ satisfies \skip 0.2in

(i) (K.1) and (K.2), if $0<\frac{2+l}{p-1}<\frac{n-2}{2}$ or

(ii) (K.1) and (K.3), if $\frac{n-2}{2}<\frac{2+l}{p-1}<n-2$.

Then,
\begin{eqnarray*}\lim_{r\to\infty}r^{m}u(r)=u_{\infty} \equiv \left\{
\begin{array}{ll}
[\frac{2+l}{p-1}(n-2-\frac{2+l}{p-1})]^{\frac{1}{p-1}}/k_{\infty}^{\frac{1}{p-1}}\ \ \mbox{or}\\\\
0.
\end{array}
\right.
\end{eqnarray*}
Furthermore, if $u_{\infty}=0$, then $\lim r^{n-2}u(r)$ exists and
is finite and positive.\vskip 0.2in

{\bf Remark 1.1.}\ When $l=-2$, a result similar to Theorem A holds
(see [16,17]).\vskip 0.1in

 The equation (1.2), if  $l_{1}=l_{2}=0$,
reduces to the following famous Lane-Emden equation
\begin{equation}
\Delta u +u^{p}+u^{q}=0 \ \ \ \ \ \mbox{in} \ \ \ R^{n},\ n\geq3.
\end{equation} \vskip0.2in

The equation (1.4) has been paid much attention recently. When $p$
and $q$ are in the range $\frac{n}{n-2}<p<\frac{n+2}{n-2}<q$, the
mixed growth structure (supercritical for $u$ large and subcritical
for $u$ small) has a profound impact on the existence and
non-existence theory, and changes the outcome for the Lane-Emden
equation, the analysis is surprisingly difficult. Recently Bamon etc
in [5] proved that if $q$ is fixed and let $p$ approach
$\frac{n+2}{n-2}$ from below, then (1.4) has a large number of
radial solutions. A similar result holds in [5] for
$p>\frac{n}{n-2}$ and $p$ is fixed,
 while letting $q$ approach $\frac{n+2}{n-2}$. In
addition, they proved that  (1.4) don't possess any solution if $q$
is fixed and then let $p$ be close enough to $\frac{n}{n-2}$.

For physical reasons, we consider the positive radial solutions of
equation (1.2), it reduces to
\begin{equation}
u''+\frac{n-1}{r}u'+r^{l_{1}}u^{p}+r^{l_{2}}u^{q}=0\ \ n\geq3.
\end{equation}
where $ r=|x|$.\vskip 0.1in

In order to state our main results concerning this question, we need
some definitions, which will be used throughout this paper:
$$
\begin{array}{lll}
\alpha_{1}=\frac{2+l_{1}}{p-1}, &
\alpha_{2}=\frac{2+l_{2}}{q-1},\\\\
\lambda_{1}^{p-1}=\alpha_{1}(n-2-\alpha_{1}),&
\lambda_{2}^{q-1}=\alpha_{2}(n-2-\alpha_{2}),\\\\
v_{1}(r)=r^{\alpha_{1}}u(r),& v_{2}(r)=r^{\alpha_{2}}u(r),\ \

\end{array}
$$
(note that $\alpha_{1}>\alpha_{2}, n-2-\alpha_{1}>0$).

We sometimes suppose $$  p\neq\frac{n+2+2l_{1}}{n-2},\ \ q
\neq\frac{n+2+2l_{2}}{n-2}. \eqno(1.6)
$$

Now we give the following definitions:\vskip0.1in

$u(r)$ is said to be singular solution of (1.2) at infinity if
$\limsup_{r\to\infty}r^{\alpha_{1}}u(r)>0$, and is said to be
regular at infinitely if $\lim_{r\to\infty} r^{n-2}u$ exists.
Similarly, $u(r)$ is said to be singular solution of (1.2) at 0 if
$\limsup_{r\to0}r^{\alpha_{2}}u>0$, is said to be regular at 0 if
$\lim _{r\to0}u$ exists.\vskip0.1in

 One of the main results is  as follows.\vskip0.2in

{\bf Theorem 1.}\ \ If $-2<l_{2}< l_{1}\leq0, \frac{n+l_{1}}{n-2}<p<
q$,
 and $p,q$ satisfy (1.6)   Then, for any solution of (1.3), we
have:

$$\lim_{r\to\infty}r^{\alpha_{1}}u=\lambda_{1},\  \mbox{or}\
\lim_{r\to\infty}r^{n-2}u=c_{1}\eqno(1.7)$$ for some constant
$c_{1}>0$. Moreover,
$$\lim_{r\to0}r^{\alpha_{2}}u=\lambda_{2},\   \mbox {or} \ \lim_{r\to 0}u=c_{2}\eqno(1.8)$$
for some constant $c_{2}>0$.\vskip0.2in

{\bf Remark 1.2.}\ If, say, the origin is a singularity, the term
$|x|^{l_{2}}u^{q}$ is dominant and the term $|x|^{l_{1}}u^{p}$ is
thus a small perturbation. At infinity, the situation is just
reversed, i.e. the term $|x|^{l_{1}}u^{p}$ is dominant and the term
$|x|^{l_{2}}u^{q}$ is a small perturbation.\vskip 0.2in

The device in the proof of Theorem 1 is an energy function which
plays a central role in the convergence of (1.7) and (1.8). If
$p=\frac{n+2+2l_{1}}{n-2}$, then the coefficient of the energy
function is 0 and $r^{\frac{2+l_{1}}{p-1}}u$ could oscillates
endlessly near the $\infty$ . As the same reason, if
$q=\frac{n+2+2l_{2}}{n-2}$, then $r^{\frac{2+l_{2}}{q-1}}u$ could
oscillates endlessly near the origin.\vskip0.1in

{\bf Theorem 2.}  If $-2<l_{2}< l_{1}\leq0, \frac{n+l_{1}}{n-2}<p<
q$, let $u$ be a positive radial solution of (1.2), then \vskip0.1in

(i) if $q=\frac{n+2+2l_{2}}{n-2}$, then either (1.8) holds, or
$v_{2}(r)$ oscillates endlessly near the origin between two
sequences ${\mu_{1,i}}$ and ${\mu_{2,i}}$ satisfying
$0<\mu_{1,i}<\mu_{2,i}$ and $$\lim_{i\to\infty}\mu_{1,i}=\mu_{1}, \
\ \ \ \ \lim_{i\to\infty}\mu_{2,i}=\mu_{2}$$
$$\mu_{1}=\lim_{r\to0}\inf r^{\frac{n-2}{2}}u(r)<\lim_{r\to0}\sup
r^{\frac{n-2}{2}}u(r)=\mu_{2},$$ where $\mu_{1}$ and $\mu_{2}$ are
fixed values satisfying
$$0<\mu_{1}\leq\lambda_{2}\leq\mu_{2}, \ \ b(\mu_{1})=b(\mu_{2})<0,\ \
\mbox{with}\ \ b(v)=\frac{1}{q+1}v^{q+1}-\frac{\lambda_{2}^{q-1}
}{2}v^{2}; \eqno(1.9)$$

(ii) if $p=\frac{n+2+2l_{1}}{n-2},$ then either (1.7) holds or
$v_{1}(r)$ oscillates endlessly near $\infty$ between two sequences
${\mu'_{1,i}}$ and ${\mu'_{2,i}}$ satisfying
$0<\mu'_{1,i}<\mu'_{2,i}$ and
$$\lim_{i\to\infty}\mu'_{1,i}=\mu'_{1}, \ \ \ \ \
\lim_{i\to\infty}\mu'_{2,i}=\mu'_{2}$$
$$\mu_{1}'=\lim_{r\to\infty}\inf r^{\frac{n-2}{2}}u(r)<\lim_{r\to\infty}\sup
r^{\frac{n-2}{2}}u(r)=\mu_{2}',$$ where $\mu_{1}'$ and $\mu_{2}'$
are fixed values satisfying
$$0<\mu_{1}'\leq\lambda_{1}\leq\mu_{2}', \ \ b_{1}(\mu_{1}')=b_{1}(\mu_{2}'),\ \
\mbox{with}\ \
b_{1}(v)=\frac{1}{p+1}v^{p+1}-\frac{\lambda_{1}^{p-1}}{2}v^{2}.
\eqno(1.10)$$\vskip0.1in

For equation (1.3), when $K\equiv1$, it becomes into the simplest
model which is the generalized Lan-Emden equation or Emden-Fowler
equation in astrophysics
$$\Delta u +u^{p}=0\ \ \ \ \ \mbox{in}\ \ R^{n},\ \ n\geq 3.
\eqno(1.11)$$\vskip0.1in

For Equation (1.11), we have the following uniqueness
result.\vskip0.1in

 {\bf Theorem B$^{[22]}$.} If $\frac{n}{n-2}<p\leq \frac{n+2}{n-2}$, then
 (1.11) admits exactly one solution
 $\frac{2}{p-1}(n-2-\frac{2}{p-1})r^{-\frac{2}{p-1}}$, singular at
 the infinity. If $p>\frac{n+2}{n-2}$ then (1.11) admits exactly one
 solution $\frac{2}{p-1}(n-2-\frac{2}{p-1})r^{-\frac{2}{p-1}}$,
 singular at the origin.\vskip0.1in

{\bf Remark 1.3.}\ From the Theorem B, we know that (1.11) has only
one solution that is singular at infinity for $\frac{n}{n-2}<p\leq
\frac{n+2}{n-2}$ and at the origin for $p>\frac{n+2}{n-2}$. Inspired
by the Theorem B, we derive the following result.\vskip0.1in

{\bf Theorem 3.} (i) if  $-2<l_{2}< l_{1}\leq0,
\frac{n+l_{1}}{n-2}<p< q<\frac{n+2+2l_{2}}{n-2}$, then (1.2) admits
exactly one solution singular at infinity. This solution has the
following exact limits:
$$\lim_{r\to\infty}r^{\alpha_{1}}u=\lambda_{1};\ \ \ \lim_{r\to
0}r^{\alpha_{2}}u=\lambda_{2},$$

(ii) if  $-2<l_{2}< l_{1}\leq0, \frac{n+2+2l_{1}}{n-2}<p< q$,
 then (1.2) admits exactly one solution
singular at the origin. This solution has the exact limits:
$$\lim_{r\to\infty}r^{\alpha_{1}}u=\lambda_{1}, \ \ \
\lim_{r\to0}r^{\alpha_{2}}u=\lambda_{2}.$$

This paper is organized as follows. In Section 2, we make some basic
observations and fundamental estimates of positive solutions. In
Section 3, the asymptotic behaviors at $\infty$ and 0 of positive
solutions are studied. Finally, the uniqueness result is
established.\vskip0.1in

 Throughout this paper, unless otherwise stated, the letter $C$
will always denote various generic constant which is independent of
$u$ and change from line to line. \vskip 0.2in

{\bf 2. Preliminaries} \setcounter{section}{2}
\setcounter{equation}{0} \vskip 0.2in  In this section we present
some preliminary results for radial solutions $u(r)$ of (1.5), where
$r=|x|$ is the radius.

First we will prove the following priori estimates, which are
inspired by the work of Ni [20].\vskip 0.2in

{\bf Lemma 2.1.} \  If $-2<l_{2}< l_{1}\leq0, \frac{n+l_{1}}{n-2}<p<
q$,  let $u(r)$ be a positive radial solution of (1.5) for
$r\in(0,\infty)$, then we have for some positive constant $C$\vskip
0.1in

(i) If $u$ tends to $\infty$ as $r\to0$, then $u(r)\leq C
r^{-\frac{2+l_{2}}{q-1}}$;

(ii)  If $u$ tends to 0 as $r \to \infty$, then $u(r)\leq C
r^{-\frac{2+l_{1}}{p-1}}$.
 \vskip0.2in

{\bf Proof.} \ For a radial solution $u=u(r)$, we rewrite (1.2) in
the
 following form:
 $$
 (r^{n-1}u')'+r^{n-1}(r^{l_{1}}u^{p}+r^{l_{2}}u^{q})=0. \eqno
 (2.1).$$
 Since $u\to \infty$ as $r\to 0$,  there exists small $r_{0}>0$, such that $u'<0$ in $(0,r_{0})$.

 Integrating (2.1) from $\bar{r}$ to $r$ ($\bar{r}<r<r_{0}$), we
 obtain
$$
r^{n-1}u'(r)=\bar{r}^{n-1}u'(\bar{r})-\int_{\bar{r}}^{r}s^{n-1}(s^{l_{1}}u^{p}+s^{l_{2}}u^{q})ds.$$
Therefore, $
r^{n-1}u'(r)\leq-\int_{\bar{r}}^{r}s^{n-1}(s^{l_{1}}u^{p}+s^{l_{2}}u^{q})ds$,
\ for\  all\ $0<\bar{r}<r$. Then letting $\bar{r}\to 0$, we obtain
\begin{eqnarray*}
r^{n-1}u'(r)&\leq &-\int_{0}^{r}s^{n-1}(s^{l_{1}}u^{p}+s^{l_{2}}u^{q})ds \\
&\leq& -\int_{0}^{r}s^{n-1}s^{l_{2}}u^{q}ds.
\end{eqnarray*}
Since $u$ is decreasing near $r=0$, we find that
$$r^{n-1}u'(r)<-u^{q}(r)\int_{0}^{r}s^{n-1}s^{l_{2}}ds=-\frac{1}{n+l_{2}}r^{n+l_{2}}u^{q}(r),$$
which in turn leads to $$\frac{u'(r)}{u^{q}(r)}\leq-
r^{l_{2}+1}.\eqno(2.2)$$ Integrating (2.2) over $(\bar{r},r)$, we
have
$$\int_{\bar{r}}^{r}\frac{u'(s)}{u^{q}}ds\leq
-\int_{\bar{r}}^{r}s^{l_{2}+1}ds.\eqno(2.3)$$ It follows from (2.3)
that
$$u^{1-q}(r)\geq
u(\bar{r})^{1-q}+\frac{q-1}{l_{2}+2}(r^{l_{2}+2}-\bar{r}^{l_{2}+2}).$$
Letting $\bar{r}\to 0$, we have
$$u^{1-q}(r)\geq Cr^{2+l_{2}}.$$ So we have $u(r)\leq
Cr^{\frac{2+l_{2}}{1-q}} \ \ \mbox{for}\  0<r<r_{0}$, and the proof
is completed. \vskip 0.2in

Part (ii) of Lemma 2.1 may be handled in a similar way. Similar to
the proof of (i) there exists a large number $R>0$ such that, for
all $r>R$,
$$r^{n-1}u'(r)<-\int_{R}^{r}s^{n-1}s^{l_{1}}u^{p}ds.$$
By a similar computation, we have $$\frac{u'(r)}{u^{p}}\leq
C(\frac{R^{n+l_{1}}}{n+l_{1}}r^{1-n}-r^{1+l_{1}}).\eqno(2.4)$$
Integrating (2.4) over $(R,r)$, we have $$ u^{1-p}(r)\geq
Cr^{2+l_{1}},\ \ \mbox{at}\ \ r=\infty.$$ So we have $u(r)\leq
Cr^{\frac{2+l_{1}}{1-p}} \ \mbox{as}\ r\to \infty$, and the proof of
part (ii) is completed.\vskip 0.2in

{\bf Lemma 2.2.}\  Let $u$ be a positive radially symmetric solution
of (1.2), then there exist two positive number $\bar{r}$ and
$\bar{\bar{r}}$ such that \vskip 0.1in

(i) $|u'(r)|\leq Cr^{-(\frac{2+l_{2}}{q-1}+1)};\ \ \ |u''(r)|\leq
Cr^{-(\frac{2+l_{2}}{q-1}+2)}$ \ \ for $0<r<\bar{r}$.\vskip 0.1in

(ii) $|u'(r)|\leq Cr^{-(\frac{2+l_{1}}{p-1}+1)};\ \ \  |u''(r)|\leq
Cr^{-(\frac{2+l_{1}}{p-1}+2)}$ \ \ for $r>\bar{\bar{r}}$.\vskip
0.2in

{\bf Proof.}\ (i)\ \  Integrate (1.2) in a small ball $B_{r}$ with
radius $r$ centered at 0. From the Lemma 2.1 and Green's identity,
we obtain
\begin{eqnarray*}
-\omega_{n}r^{n-1}u'(r)&=&-\int_{B_{r}}\Delta u
=\int_{B_{r}}(|x|^{l_{1}}u^{p}+|x|^{l_{2}}u^{q})dx\\
&\leq& C \int_{0}^{r}s^{l_{2}}u^{q} s^{n-1}ds\leq
C\int_{0}^{r}s^{-\frac{2+l_{2}}{q-1}q+l_{2}+n-1}ds\\
&=&Cr^{-\frac{2+l_{2}}{q-1}q+l_{2}+n}.
\end{eqnarray*}
So as $r\to 0$, we have $$-u'(r)=|u'(r)|\leq
Cr^{-\frac{2+l_{2}}{q-1}q+l_{2}+1}=Cr^{-(\frac{2+l_{2}}{q-1}+1)}$$
\begin{eqnarray*}
|u''|&\leq&\frac{n-1}{r}|u'| +r^{l_{1}}u^{p}+r^{l_{2}}u^{q}\\
&\leq&C[r^{-(\frac{2+l_{2}}{q-1}+2)}+r^{l_{2}-\frac{(2+l_{2})q}{q-1}}]\leq
Cr^{-(\frac{2+l_{2}}{q-1}+2)}.
\end{eqnarray*}
And the proof of (i) is completed.\vskip0.1in

Part (ii) of Lemma 2.2 may be handled in a similar fashion. As (i),
for large $r$ we have
\begin{eqnarray*}
\omega_{n}r^{n-1}(-u'(r))&=&-\int_{B_{r}}\Delta u
=\int_{B_{r}}s^{l_{1}}u^{p}+s^{l_{2}}u^{q}dx\\
&\leq& C+C\int_{R_{0}}^{r}s^{-\frac{2+l_{1}}{p-1}p+l_{1}+n-1}ds,
\end{eqnarray*}
where $R_{0}$ is a large positive number. By a simple computation,
we obtain $$|u'(r)|\leq Cr^{-(\frac{2+l_{1}}{p-1}+1)},\ \ \
|u''(r)|\leq Cr^{-(\frac{2+l_{1}}{p-1}+2)}\ \ \  \mbox{at}\ \
r=\infty,$$ and the proof is over.\vskip0.1in

 {\bf Lemma 2.3.} Suppose that $u$ is a positive solution of (1.5).
 Let $v(r)=r^{\alpha}u(r)$, then $v$ satisfies
 $$v''+\frac{n-1-2\alpha}{r}v'-\frac{(n-2-\alpha)\alpha}{r^{2}}v+r^{l_{1}-(p-1)\alpha}v^{p}+r^{l_{2}-(q-1)\alpha}v^{q}=0.$$
Let $\alpha=\alpha_{1}$, then we have ($v_{1}=r^{\alpha_{1}}u$)
$$v_{1}''+\frac{n-1-2\alpha_{1}}{r}v_{1}'-\frac{(n-2-\alpha_{1})\alpha_{1}}{r^{2}}v_{1}+\frac{v_{1}^{p}}{r^{2}}+r^{l_{2}-(q-1)\alpha_{1}}v_{1}^{q}=0.\eqno(2.5)$$
Let $\alpha=\alpha_{2}$, then we have ($v_{2}=r^{\alpha_{2}}u$)
$$v_{2}''+\frac{n-1-2\alpha_{2}}{r}v_{2}'-\frac{(n-2-\alpha_{2})\alpha_{2}}{r^{2}}v_{2}+\frac{v_{2}^{q}}{r^{2}}++r^{l_{1}-(p-1)\alpha_{2}}v_{2}^{p}=0\eqno(2.6)$$

This lemma can be proved by straight forward calculations, thus we
omit it here.\vskip0.1in

 {\bf Lemma 2.4.} we have $$v_{1}'^{2}r\in L^{1}(R,\infty), \ \ \ v_{2}'^{2}r\in L^{1}(0,R)
 ,\eqno(2.7)$$ where $R$ is a large positive number.\vskip0.1in

{\bf Proof.} Multiplying (2.5) by $v_{1}'r^{2}$ and integrating from
$R$ to $r>R$, we obtain
$$ \frac{v_{1}'^{2}s^{2}}{2}\Bigg|_{R}^{r}+c^{1}\int_{R}^{r}v_{1}'^{2}sds-\frac{\lambda_{1}^{p-1}}{2}v_{1}^{2}\Bigg|_{R}^{r}
+\frac{1}{p+1}v_{1}^{p+1}\Bigg|_{R}^{r}+\int_{R}^{r}s^{l_{2}-\alpha_{1}(q-1)+2}v_{1}'v_{1}^{q}ds=0,\eqno
(2.8)$$ where $c^{1}=n-2-2\alpha_{1}.$

 From (ii) of Lemma 2.2, we obtain
$$|v_{1}'|\leq\frac{C}{r};\ \ |v_{1}''|\leq\frac{C}{r^{2}},\ \ \mbox{at}\ \
r=\infty,\eqno(2.9)$$ so $ v_{1}^{p+1}\Bigg|_{R}^{r},
\frac{v_{1}'^{2}s^{2}}{2}\Bigg|_{R}^{r}$ and \
$\int_{R}^{r}s^{l_{2}-\alpha_{1}(q-1)+2}v_{1}'v_{1}^{q}ds$ are
bounded at $r=\infty$, and from that we have
$$ \int_{R}^{r}v_{1}'^{2}sds\leq C,\ \ \ \mbox{for all}\ \ r>R$$
since $c^{1}\neq 0$ by (1.6), and $(2.7)_{1}$ follows. $(2.7)_{2}$
is handled by the similar way, we omit it here.

{\bf Lemma 2.5.}  We have
$$\lim_{r\to\infty}rv_{1}'=0;\ \
\lim_{r\to0}rv_{2}'=0.\eqno(2.10)$$

{\bf Proof.} \ Now we prove the $(2.10)_{1}$. Suppose for
contradiction that it is not true, then there exist a sequence
${r_{k}}\to +\infty$  such that
$$|v_{1}'(r_{k})r_{k}|\geq C.$$
From (2.9), one obviously has, near $\infty$,
$$|(v_{1}'^{2}r^{2})'|\leq\frac{M}{r},$$ for some $M>0$. Combining
the above two inequalities yields
$$|v_{1}'^{2}(r)r^{2}-v_{1}'(r_{k})^{2}r_{k}^{2}|\leq
M|r-r_{k}|\max(\frac{1}{r},\frac{1}{r_{k}}),$$ and so
$$v_{1}'^{2}(r)r^{2}\geq\frac{C^{2}}{2},\ r\in
[(1+\varepsilon)^{-1}r_{k},(1+\varepsilon)r_{k}], \ \
\varepsilon(1+\varepsilon)=\frac{C^{2}}{2M}.$$  This contradicts
$(2.7)_{1}$.

The proof of $(2.10)_{2}$ is handled by the same way, as the
previous proof, we have, near 0,
$$|(v_{2}'^{2}r^{2})'|\leq\frac{M}{r}$$ for some $M>0$, and  by a
similar caculation, we will obtain a contradiction and the proof is
completed.\vskip0.1in

{\bf Lemma 2.6} \ Let $u$ be a positive superharmonic function near
$\infty$ and $\bar{u}$ its spherical mean. Then, $r^{n-2}\bar{u}$ is
increasing as $r\to \infty$.\vskip0.1in

{\bf Proof.} Put $f(t):=r^{n-2}\bar{u}(r), t=\log r.$ Then $f$
satisfies
$$f''-(n-2)f'\leq0$$
and $f'(t)\leq e^{(n-2)(t-T)}f'(T)$ on [T,t] for $T$ large. Because
 $f$ is positive, $f$ must be increasing near $\infty$. It implies
that $(r^{n-2}\bar{u}(r))_{r}>0$ near $\infty$.\vskip0.2in

{\bf 3. Asymptotic behavior} \setcounter{section}{3}
\setcounter{equation}{0}\vskip0.1in

 \ \ In this section we investigate the asymptotic behavior at
 $\infty$ and  0 of positive radial solutions of (1.2). We will
 prove any positive radial solution of (1.2) must behave either like
 $r^{-\alpha_{1}}$ or $r^{2-n}$ at $\infty$. In addition, if any positive radial
  solution of (1.2) is singular at the origin, then it must behave
  like $r^{-\alpha_{2}}$. Now we give the proof of Theorem
  1.\vskip0.1in

 {\bf The proof of Theorem 1.} Consider the function
$a(r)=\frac{v_{1}^{p+1}}{p+1}-\frac{\lambda_{1}^{p-1}v_{1}^{2}}{2}.$
By Lemma 2.4 and Lemma 2.5, we have, for fixed $R>0$,
$$ \frac{v_{1}'^{2}r^{2}}{2}\to 0,\ \int_{R}^{r}v_{1}'^{2}s \to
c_{2},\
\int_{R}^{r}r^{l_{2}-\frac{(2+l_{1})(q-1)}{p-1}+2}v_{1}'v_{1}^{p}\to
c_{3},\ \ \mbox{as}\ r\to \infty$$ for some constant $c_{2}$ and
$c_{3}$. This implies, by (2.8), that $a(r)$ must tend to a finite
constant $c_{4}$ as $t\to \infty$. We claim $v_{1}$ approaches  a
finite limit as $r\to \infty$. If not, we may choose two sequences
$\{\eta_{i}\}$ and $\{\xi_{i}\}$ going to $\infty$ as $i\to\infty$
such that\vskip0.1in
\begin{eqnarray*}
\left\{
\begin{array}{ll}
\{\eta_{i}\}\ \mbox{are local minima of}\ v_{1}, \{\xi_{i}\} \
\mbox{are local
maxima}.\\
\eta_{i}<\xi_{i}<\eta_{i+1}, i=1,2, ....
\end{array}
\right.
\end{eqnarray*}
And we have
\[v_{1}(\eta_{i})\to m_{1},\ \ v_{1}(\xi_{i})\to m_{2}\ \ \
\mbox{as}\ i \to \infty, \] for some positive constants $m_{1},
m_{2}$ with $m_{1}< m_{2}$. Since $a(r)$ tend to a finite constant
as $ r\to \infty$, we have $$
\frac{m_{1}^{p+1}}{p+1}-\frac{\lambda_{1}^{p-1}m_{1}^{2}}{2}=
\frac{m_{2}^{p+1}}{p+1}-\frac{\lambda_{1}^{p-1}m_{2}^{2}}{2}=c_{4},$$
where $c_{4}$ is a constant. However, the intermediate value theorem
shows that there exists $r_{i}\in (\eta_{i},\xi_{i})$ such that
$$v(r_{i})=m_{0},\ \  m_{1}<m_{0}<m_{2}\ \mbox{and}\  \frac{da(v)}{dv}(r_{i})=0.$$
Furthermore, we have
$\frac{m_{0}^{p+1}}{p+1}-\frac{\lambda_{1}^{p-1}m_{0}^{2}}{2}\neq
c_{4}$, since
$\frac{v_{1}^{p+1}}{p+1}-\frac{\lambda_{1}^{p-1}v_{1}^{2}}{2}$ has
only one minima for $v_{1}\in[0,+\infty]$.
  A contradiction is obtained. Similarly, we conclude that
$r^{\alpha_{2}}u(r)$ approaches a finite limit as $r\to
0$.\vskip0.1in

Claim. (1) $\lim_{r\to\infty}v_{1}(r)$ must be either 0 or
$\lambda_{1}$;

\ \ \ \ \ \ \ \ \ \  (2) $\lim_{r\to0}v_{2}(r)$ must be either 0 or
$\lambda_{2}$.\vskip0.1in

We only prove (1), the demonstration of (2) being the same. We
denote the $\lim_{r\to\infty}v_{1}(r)$ by $v_{\infty}$. Now, if
$v_{\infty}\neq 0$, we want to show that $v_{\infty}=\lambda_{1}$.
From (2.5), we have
$$\frac{d^{2} v_{1}}{d t^{2}}+(n-2-2\alpha_{1})\frac{d v_{1}}{d t}
-(n-2-\alpha_{1})\alpha_{1}v_{1}+v_{1}^{p}+e^{(l_{2}-(q-1)\alpha_{1}+2)t}v_{1}^{q}=0,$$
where $t=\log r$. From Lemma 2.5, we have $v'(t)\to0$ as
$t\to\infty$. So $\lim_{t\to\infty} v''(t)$  exists and must be 0.
Immediately, we have $v_{\infty}=\lambda_{1}$ or 0.

If $\lim_{r\rightarrow\infty}v_{1}(r)=0$ or
$\lim_{r\rightarrow0}v_{2}(r)=0$, (2.5) and (2.6) suggests that
$v_{i} (i=1,2)$ tends to zero at an algebraic rate. Indeed, since
$p,q>1$ and $v_{i} (i=1,2)$ is expected to satisfy asymptotically
the following equations\vskip0.1in
\begin{eqnarray*}
\begin{array}{ll}
v_{1}''+\frac{n-1-2\alpha_{1}}{r}v_{1}'-\frac{\lambda_{1}^{p-1}}{r^{2}}v_{1}=0
& \mbox{as}\ \ r \to \infty, \\\\
v_{2}''+\frac{n-1-2\alpha_{2}}{r}v_{2}'-\frac{\lambda_{2}^{q-1}}{r^{2}}v_{2}=0
& \mbox{as}\ \ r \to 0.
\end{array}
\end{eqnarray*}
Therefore $v_{i}$ should satisfy the following asymptotical
behaviors
\[v_{1}\approx r^{-(n-2-\alpha_{1})}\ \mbox{at}\ \infty,\ \ \
v_{2}\approx r^{\alpha_{2}}\  \mbox{at}\ 0. \]

Now we claim: (i) If $\lim_{r\to\infty}v_{1}=0$, then
$\lim_{r\to\infty}r^{n-2}u=c>0$;

\ \ \ \ \ \ \ \ \ \ \ \ \ \ \ \ \ \ \ \ (ii) If $\lim_{r\to 0}
v_{2}=0$, then $\lim_{r\to 0}u=c_{1}>0$.\vskip0.1in

First, we prove (i), by assumption, for any $\varepsilon>0$ there
exists a positive number $r_{\varepsilon}$ such that $v_{1}$
satisfies
$$v_{1}''+\frac{n-1-2\alpha_{1}}{r}v_{1}'-\frac{(\lambda_{1}^{p-1}-\varepsilon)}{r^{2}}v_{1}\geq
0, \ r>r_{\varepsilon}. \eqno(3.11)$$ The characteristic equation of
(3.11) has the two characteristic values
\begin{eqnarray*}
\begin{array}{ll}
a_{1}=\alpha_{1}-\frac{n-2-\sqrt{(n-2)^{2}-4\varepsilon}}{2}=\alpha_{1}+O(\varepsilon)\\
a_{2}=\alpha_{1}-\frac{n-2+\sqrt{(n-2)^{2}-4\varepsilon}}{2}=\alpha_{1}+2-n+O(\varepsilon)
\end{array}
\end{eqnarray*}
Rewrite (3.11) $$(D-\frac{a_{1}-1}{r})(D-\frac{a_{2}}{r})v_{1}\geq
0,$$  where  $D:=\frac{d}{dr}, D^{2}:=\frac{d^{2}}{dr^{2}}.$

Let $(D-\frac{a_{2}}{r}v_{1})=U_{1}$, so we have
$$U_{1}'+\frac{1-a_{1}}{r}U_{1}\geq 0,$$
from which, we have
$$[r^{1-\alpha_{1}+O(\varepsilon)}(D-\frac{a_{2}}{r})v_{1}]'\geq 0.\eqno(3.12)$$
Observe that, for $\varepsilon$ small enough,
$$\lim_{r\to\infty}r^{1-\alpha_{1}+O(\varepsilon)}(D-\frac{a_{2}}{r})v_{1}=0
\ \ \ \mbox{by}\ \ (2.9),$$ since $1-\alpha_{1}+O(\varepsilon)<1$.
It follows from (3.12) that
$$(D-\frac{a_{2}}{r})v_{1}\leq 0, \ \
r>r_{\varepsilon}.$$ Integrating once from $r_{\varepsilon}$ to $r$
yields
$$v_{1}\leq
c_{\varepsilon}r^{a_{2}}=c_{\varepsilon}r^{\alpha_{1}+2-n+O(\varepsilon)},\
\ \  r>r_{\varepsilon}.$$ Thus for $\varepsilon$ sufficiently small,
$$v_{1}''+\frac{n-1-2\alpha_{1}}{r}v_{1}'-\frac{\lambda_{1}^{p-1}}{r^{2}}v_{1}=g(r)\eqno(3.13)$$
with
$$g(r)=\frac{v_{1}^{p}}{r^{2}}+r^{l_{2}-(q-1)\alpha_{1}}v_{1}^{q}=
O(r^{-2-\delta})\ \  \mbox{near}\ \infty,\eqno(3.14)$$ where
$\delta$ is a positive constant.

 Applying the method of variation of parameters to (3.13), $v_{1}$ is represented
by
\begin{eqnarray*}
v_{1}(r)&=&C_{1}(R)r^{\alpha_{1}}+C_{2}(R)r^{\alpha_{1}+2-n}\\\\
&+&\frac{r^{\alpha_{1}+2-n}}{2-n}\int_{R}^{r}s^{n-1-\alpha_{1}}g(s)ds-
\frac{r^{\alpha_{1}}}{2-n}\int_{R}^{r}s^{1-\alpha_{1}}g(s)ds.
\end{eqnarray*}
By a similar computation, we have, from (3.14) and the bounded of
$v_{1}(r)$ ,
$$v_{1}(r)=Cr^{\alpha_{1}+2-n}+o(r^{\alpha_{1}+2-n}),$$
hence,
$$r^{n-2}u(r)\leq C \ \ \ \mbox{at}\ \ r=\infty.$$
By Lemma 2.6, there exists a constant $c$ such that $r^{n-2}u(r)\to
c$ as $r\to \infty$.

The proof of (ii) is handled by the same way. As the proof of (i),
for any $\varepsilon'>0$ there exists a positive number
$r_{\varepsilon'}$ such that
$$\lim_{r\to0}r^{n-1-\alpha_{2}+O(\varepsilon')}(D-\frac{\alpha_{2}}{r})v_{2}=0,$$
and
$$(D-\frac{\alpha_{2}}{r})v_{2}\geq0\ \ \ 0<r<r_{\varepsilon'}.$$
Integrating once from $r$ to $r_{\varepsilon'}$ yields
$$v_{2}\leq c_{\varepsilon'}r^{\alpha_{2}+O(\varepsilon')}.$$
Applying the method of variation of parameters to (2.6), we
immediately that $v_{2}$ is bounded by $r^{\alpha_{2}}$, and in turn
$u$ is bounded, by standard theory, we have $\lim_{r\to
0}u(r)=c_{1}$. \vskip0.1in

When $p=\frac{n+2+2l_{1}}{n-2}$ or $q=\frac{n+2+2l_{2}}{n-2}$, there
is another possibility; $r^{\frac{n-2}{2}}u$ could oscillate
endlessly  near $\infty$ or oscillate endlessly near the
origin.\vskip0.1in

{\bf The proof of Theorem 2.} (i) The transformed function
$v_{2}(t):=r^{\alpha_{2}}u, t= \log r,$ satisfies
$$v_{2}''-\lambda_{2}^{q-1}v_{2}+v_{2}^{q}+e^{\delta_{1}t}v_{2}^{p}=0,$$
where $\delta_{1}=(\alpha_{1}-\alpha_{2})(p-1)>0.$ Note that by
Lemma 2.1 (i), $v_{2}(t)$ is bounded near $-\infty$. Define an
energy function
$$E(t):=
\frac{1}{2}v_{2}'^{2}-\frac{\lambda_{2}^{q-1}}{2}v_{2}^{2}+\frac{1}{q+1}v^{q+1}.$$
As in (2.8), we have
$$E(t)=C(T)+\int_{t}^{T}e^{\delta_{1}s}v_{2}^{p}v_{2}'ds,$$ where
$T$ is a fixed  number.

In case that $v_{2}$ oscillates near $-\infty$, then we may suppose
that
$$0\leq\mu_{1}=\lim_{t\to-\infty}\inf
v_{2}(t)<\lim_{t\to-\infty}\sup v_{2}(t)=\mu_{2}<\infty.$$ Then,
there exists two sequences ${\eta_{i}}$ and ${\varepsilon_{i}}$
going to $-\infty$ as $i\to\infty$ such that ${\eta_{i}}$ and
${\varepsilon_{i}}$ are local minima and local maxima of $v_{2}$,
respectively, satisfying $\eta_{i}<\varepsilon_{i}<\eta_{i+1},
i=1,2,....$ From (2.8), we know
$\lim_{t\to-\infty}\int_{t}^{T}e^{\delta_{1}s}v_{2}^{p}v_{2}'ds$
exists, so $E=\lim_{t\to-\infty}E(t)$ exists, and from that we have
$$\lim_{i\to\infty}E(\eta_{i})=-\frac{\lambda_{2}^{q-1}}{2}\mu_{1}^{2}+\frac{1}{q+1}\mu_{1}^{q+1}
=-\frac{\lambda_{2}^{q-1}}{2}\mu_{2}^{2}+\frac{1}{q+1}\mu_{2}^{q+1}=\lim_{i\to\infty}E(\varepsilon_{i}),$$
which implies $b(\mu_{1})=b(\mu_{2})$.

Observe that for each $i>1$,
$$0\leq
v_{2}''(\eta_{i})=\lambda_{2}^{q-1}v_{2}(\eta_{i})-v_{2}^{q}(\eta_{i})-v_{2}^{p}(\eta_{i})e^{\delta_{1}\eta_{i}},\eqno(3.16)$$
while
$$0\geq
v_{2}''(\varepsilon_{i})=\lambda_{2}^{q-1}v_{2}(\varepsilon_{i})-
v_{2}^{q}(\varepsilon_{i})-v_{2}^{p}(\varepsilon_{i})e^{\delta_{1}\varepsilon_{i}}.\eqno(3.17)$$
Putting $v_{2}(\eta_{i})=\mu_{1,i},
v_{2}(\varepsilon_{i})=\mu_{2,i}$, we have
$$\lim_{i\to\infty}\mu_{1,i}=\mu_{1}, \ \ \ \ \
\lim_{i\to\infty}\mu_{2,i}=\mu_{2}.$$ From (3.16) and (3.17) we
obtain
$$0<\mu_{1}\leq\lambda_{2}\leq\mu_{2}.$$

The proof of (ii) is handled by the similar way, we omit it
here.\vskip0.2in

{\bf 4. A uniqueness result} \setcounter{section}{4}
\setcounter{equation}{0}\vskip0.1in

In this section, we shall prove a uniqueness result for positive
radial singular solutions of (1.2) when $p>\frac{n+2+2l_{1}}{n-2}$
and $q<\frac{n+2+2l_{2}}{n-2}$. More precisely, we require that
solutions be singular at infinity for $q<\frac{n+2+2l_{2}}{n-2}$,
and at the origin for $p>\frac{n+2+2l_{1}}{n-2}$. In order to prove
the Theorem 3, we need a finer asymptotic behavior of solution of
(1.2) near 0 and $\infty.$

Let $w(t)=v_{1}-\lambda_{1}, t=\log r$.  If
$\lim_{t\to+\infty}v(t)=\lambda_{1}$, then it satisfy $w(t)\to 0 \ \
\ \mbox{as}\ \ t\to+\infty$ and
$$w''+(n-2-2\alpha_{1})w'+(2+l_{1})(n-2-\alpha_{1})w+f(t)+(w+\lambda_{1})^{q}e^{\delta
t}=0,\eqno(4.1)$$ where
$\delta=\frac{(2+l_{1})(1-q)}{p-1}+2+l_{2}<0$, and
\begin{eqnarray*}
f(w)&=&(\lambda_{1}+w)^{p}-\lambda_{1}^{p}-p\lambda_{1}^{p-1}w=
\lambda_{1}^{p}\sum_{k=2}^{\infty}\frac{(p-k+1)}{k!}(\frac{w}{\lambda_{1}})^{k}\\
&=&\frac{(p-1)\lambda_{1}^{p-2}}{2}w^{2}+o(w^{2})\ \ \ \mbox{for}\ \
w\ \ \mbox{near} \ \ 0.  \ \ \ \ \ \ \ \ \ \ \ \ \ \ \ \ \ \ \ \ \ \
\ \ \ \
\end{eqnarray*}

{\bf Theorem 4.1.}\ \ Let  $u(r)$ be a singular solution of (1.2) at
infinity and $\frac{n+l_{1}}{n-2}<p<\frac{n+2+2l_{1}}{n-2}$, then
for any $\varepsilon\in(0,-\delta)$ we have
$$w(t)=O(e^{-\varepsilon t}),\ \ \ w'(t)=O(e^{-\varepsilon t}), \ \ \ \ \mbox{as}\ t\to+\infty.\eqno(4.2)$$

{\bf Proof.}\ For $T>0$ and $t\in (0,T)$, multiply (4.1) by $2w'(t)$
and integrate from $t$ to $T$ to obtain
\begin{eqnarray*}
&&[w'^{2}+(2+l_{1})(n-2-\alpha_{1})w^{2}]\Big|_{t}^{T}+2(n-2-2\alpha_{1})\int_{t}^{T}w'^{2}ds\\
&&+ 2\int_{t}^{T}f(t)w'+2\int_{t}^{T}(w+\lambda_{1})^{q}e^{\delta
t}w'ds=0.\ \ \ \ \ \ \ \ \ \ \ \ \ \ \ \ \ \ \ \ \ \ \ \ \ \ \ \ \ \
\ \ \ \ \ \ \ \ \ \ \ \ \ \ \ \ \ (4.3)
\end{eqnarray*}
By (2.10), it is easy to see that for large $t$
$$w'^{2}(T)\to0,\ w^{2}(T)\to0,\ \
\int_{t}^{T}f(s)w'ds=-\frac{(p-1)\lambda_{1}^{p-2}}{6}w^{3}(t)+o(w^{3}(t))\
\ \mbox{as}\ T\to+\infty.$$ It follows from Lemma 2.4 and  by
letting $T\to \infty$ in (4.3) that for large $t$
$$w'^{2}(t)+(2+l_{1})(n-2-\alpha_{1})w^{2}(t)\leq
C|w|^{3}+\int_{t}^{\infty}(w+\lambda_{1})^{q}e^{\delta s}w'ds,$$
since $n-2-2\alpha_{1}<0$. It follows that for large $t$
$$w'^{2}(t)+(2+l_{1})(n-2-\alpha_{1})w^{2}(t)\leq
C\int_{t}^{\infty}(w+\lambda_{1})^{q}e^{\delta s}w'ds\leq Ce^{\delta
t}.\eqno(4.4)$$ Hence by (4.4)
$$|w'(t)|+|w(t)|\leq Ce^{\frac{\delta}{2}t}.$$
Therefore
$$\Big|\int_{t}^{T}(w+\lambda_{1})^{q}e^{\delta t}w'ds\Big|\leq
C\int_{t}^{T}e^{\frac{3\delta}{2}s}\leq C e^{\frac{3\delta}{2}t},$$
and so by (4.4)
$$|w'(t)|+|w(t)|\leq Ce^{\frac{3\delta}{4}t}.$$
Thus for any $m>0$, using a simple iteration of $m$-step in (4.4)
yields
$$|w'(t)|+|w(t)|\leq Ce^{\frac{(2^{m}-1)\delta}{2^{m}}t},$$
and (4.2) follows by taking $m$ large.\vskip0.2in

{\bf Remark 4.1.}\  As a matter of fact, apply by the method of
variation of parameters to (4.1), the $\varepsilon$ can reach at
$-\delta$. Since from (4.2), we have $f(t)=O(w^{2}(t))=O(e^{\delta
t})$ and
$$w(t)=-\frac{1}{\omega}\int_{t}^{\infty}e^{\frac{n-2-2\alpha_{1}}{2}(s-t)}\sin\omega(s-t)[f(s)+e^{\delta
t}(v_{1}+\lambda_{1})^{p}]ds=O(e^{\delta t}),$$ where
$\omega=((2+l_{1})(n-2-\alpha_{1})-\frac{1}{4}(n-2-2\alpha_{1})^{2})^{\frac{1}{2}}.$\vskip0.2in

{\bf Remark 4.2.}\ If $u(r)$ is a singular solution of (1.2) at 0,
and $\frac{n+l_{1}}{n-2}<q<\frac{n+2+2l_{2}}{n-2}$, then we have a
similar  result
$$v_{2}(r)-\lambda_{2}=O(r^{\epsilon}),\ \
(v_{2}(r)-\lambda_{2})'=O(r^{\epsilon+1})\ \ \mbox{as}\ \ r \to 0,$$
where
$\epsilon\in(0,\frac{(2+l_{2})(1-p)}{q-1}+2+l_{1})$.\vskip0.1in

 Now we are ready to prove  Theorem 3.  Let  $u_{1}(r)$ and
$u_{2}(r)$ be two different singular solutions of (1.2) at infinity.
We introduce the function
$$\bar{w}(t)=r^{\alpha_{1}}u_{1}(r)-r^{\alpha_{1}}u_{2}(r)=w_{1}(t)-w_{2}(t),\ \ t=\log r,$$
and show that $\bar{w}(t)$ is identically zero.\vskip0.2in

{\bf Proof of Theorem 3.}\ Clearly $\bar{w}(t)$ satisfies
$$
\begin {array}{llll}
&&\bar{w}''(t)+(n-2-2\alpha_{1})\bar{w}'(t)+(2+l_{1})(n-2-\alpha_{1})\bar{w}\\
&&+(f(w_{1})-f(w_{2}))+[(w_{1}(t)+
\lambda_{1})^{q}-(w_{2}(t)+\lambda_{1})^{q}]e^{\delta t}=0.
\end{array}\eqno(4.5)
$$
For $T>0$ and $t\in(0,T),$ multiply (4.5) $2 \bar{w}'$ and integrate
from $t$ to $T$ to obtain
$$
\begin {array}{llll}
&&[\bar{w}'^{2}+(2+l_{1})(n-2-\alpha_{1})\bar{w}^{2}]\Big|_{t}^{T}+(n-2-2\alpha_{1})\int_{t}^{T}\bar{w}'^{2}\\
&&+2\int_{t}^{T}[f(w_{1})-f(w_{2})]\bar{w}'+
2\int_{t}^{T}[(w_{1}(t)+
\lambda_{1})^{q}-(w_{2}(t)+\lambda_{1})^{q}]e^{\delta s}\bar{w}'=0.
\end{array}
\eqno(4.6)
$$
Thus we have the following estimates:
$$\Big|\int_{t}^{T}[f(w_{1})-f(w_{2})]\bar{w}'\Big|\leq
C\int_{t}^{T}|\bar{w}'\bar{w}|e^{-\varepsilon s}\leq
C\int_{t}^{T}e^{-\varepsilon s}(\bar{w}^{2}+\bar{w}'^{2}),$$ and
$$\int_{t}^{T}[(w_{1}(t)+
\lambda_{1})^{q}-(w_{2}(t)+\lambda_{1})^{q}]e^{\delta s}\bar{w}'\leq
C\int_{t}^{T}e^{-\delta s}(\bar{w}^{2}+\bar{w}'^{2}).$$ For large
$t$, it follows, by letting $T\to+\infty$ and using Theorem 4.1,
that
$$\bar{w}'^{2}(t)+(2+l_{1})(n-2-\alpha_{1})\bar{w}^{2}(t)\leq
C\int_{t}^{\infty}e^{-\varepsilon s}(\bar{w}^{2}+\bar{w}'^{2}),$$
since $n-2-2\alpha_{1}\leq 0$. From which, we have
$$\bar{w}'^{2}(t)+\bar{w}^{2}(t)\leq
C\int_{t}^{\infty}e^{-\varepsilon s}(\bar{w}^{2}+\bar{w}'^{2}),$$
sine $(2+l_{1})(n-2-\alpha_{1})>0$. Hence
$\bar{w}'^{2}(t)+\bar{w}^{2}(t)\equiv 0$ for all sufficiently large
$t$ by the Gronwall inequality. According to the uniqueness of
ordinary differential equation, $w$ is identically zero for all $t$.
\vskip 0.1in

When $p>\frac{n+2+2l_{1}}{n-2}$, we use the transformation
$$\bar{w_{1}}(t)=r^{\alpha_{2}}u_{1}(r)-r^{\alpha_{2}}u_{2}(r)=z_{1}(t)-z_{2}(t),\
\ \ t=\log r,$$ where $z(t)=r^{\alpha_{2}}u(r)-\lambda_{2}\to 0$ as
$t\to-\infty$. Then $\bar{w_{1}}(t)$ satisfies
\begin{eqnarray*}
&&\bar{w_{1}}''(t)+(n-2-2\alpha_{2})\bar{w_{1}}'(t)+(2+l_{2})(n-2-\alpha_{2})\bar{w_{1}}\\
&&+(\bar{f}(z_{1})-\bar{f}(z_{2}))+[(z_{1}(t)+\lambda_{2})^{p}-(z_{2}(t)+\lambda_{2})^{p}]e^{\delta_{2}t}=0\
\ \ \ \ \ \ \ \ \ \ \ \ \ \ \ \ \ \ \ \ \ \ \ \ \ \ \ \ \ \ \ \ \
(4.7)
\end{eqnarray*}
where $\delta_{2}=(P-1)(\alpha_{1}-\alpha_{2})>0,$ and
\begin{eqnarray*}
\bar{f}(z)&=&(\lambda_{2}+z)^{q}-\lambda_{2}^{q}-q\lambda_{2}^{q-1}z=
\lambda_{2}^{q}\sum_{k=2}^{\infty}\frac{q-k+1}{k!}(\frac{z}{\lambda_{2}})^{k}\\
&=&\frac{(q-1)\lambda_{2}^{q-2}}{2}z^{2}+o(z^{2})  \ \ \mbox{for}\ t
\ \mbox{near}\ -\infty.
\end{eqnarray*}
The equation (4.7) is the same as (4.5). On the other hand, one
easily sees from previous proof that the key ingredient of the proof
for the case $p>\frac{n+2+2l_{1}}{n-2}$ is that the coefficient of
the term $\bar{w_{1}}'(t)$ is not less than 0. Fortunately, we have
$n-2-2\alpha_{2}>0$ when $p>\frac{n+2+2l_{1}}{n-2}$. Hence the proof
above carries over immediately.\vskip 0.1in

\textbf{Acknowledgements}

The authors are greatly indebted to the advisor Professor Yi Li  for
many useful discussions, suggestions and comments.

  \textit{ \\
 Baishun Lai\\
 School of Mathematics, Henan University,\\
 Kaifeng, Henan 475004, PR of China\\
 E-mail address:
 laibaishun@henu.edu.cn.\\
 Shuqing Zhou\\
School of Mathematics, Hunan Normal University\\
Changsha,  Hunan 410081, China\\
 E-mail address: zhoushuqing63@163.com\\
Qing Luo\\
School of Mathematics, Henan University,\\
 Kaifeng, Henan 475004, PR of China\\
 E-mail address:
 Lq@henu.edu.cn
}
\end{document}